\documentclass[requno]{amsart}
\usepackage[latin1]{inputenc}
\usepackage[TS1,T1]{fontenc}
\usepackage{amsfonts,amssymb,latexsym,amsmath,enumerate,amsthm}

\newtheorem{thm}{Theorem}
\newtheorem{cor}{Corollary}

\newcommand{\R}{\mathbb{R}} 
\newcommand{\p}{\partial}

\begin{document}

\title[Best Sobolev constant for the embedding $BV(\Omega)$ into $L^1(\p\Omega)$]
{Estimates of the best Sobolev constant of the embedding of $BV(\Omega)$ into $L^1(\p\Omega)$ and related shape optimization problems}

\author[N. Saintier]{Nicolas Saintier}

\address{Departamento de Matem\'atica, FCEyN UBA (1428)
\hfill\break\indent Buenos Aires, Argentina. }

\email{{\tt nsaintie@dm.uba.ar} }



\begin{abstract}
In this paper we find estimates for the optimal constant in the critical Sobolev trace inequality 
$\lambda_1(\Omega)\|u\|_{L^1(\partial\Omega)} \le \|u\|_{W^{1,1}(\Omega)}$ that are independent of $\Omega$. This estimates generalize those of \cite{BS} concerning the $p$-Laplacian to the case $p=1$. 

We apply our results to prove existence of an extremal for this embedding. We then study an optimal design problem related to $\lambda_1$, and eventually compute the shape derivative of the functional $\Omega\to\lambda_1(\Omega)$. As a consequence, we obtain that a ball of $\R^n$ of radius $n$ is critical for volume-preserving deformations.
\end{abstract}

\maketitle

Let $\Omega$ be a bounded smooth domain of $\R^N$. It is well-known that the trace embedding from $W^{1,1}(\Omega)$ into $L^1(\p\Omega)$ is continuous, where $W^{1,1}(\Omega)$ is the usual Sobolev spaces of functions $u\in L^1(\Omega$ such that $\nabla u\in L^1(\Omega)$. The best constant for this embedding is then defined by 
\begin{equation}\label{BestConst_1}
  \lambda_1(\Omega)=\inf_{u\in W^{1,1}(\Omega)\setminus W_0^{1,1}(\Omega)}
 \frac{\displaystyle \int_{\Omega} |\nabla u|\,dx + \displaystyle \int_{\Omega}|u|\, dx}{\displaystyle \int_{\partial\Omega} |u|\,dH^{N-1}}, 
\end{equation}
where $W_0^{1,1}(\Omega)$ denotes the closure for the $W^{1,1}$-norm of the space of smooth functions with compact support in $\Omega$, and $H^{N-1}$ is the $(N-1)$-dimensional Hausdorff measure. 
The purpose of this paper is to obtain estimates of $\lambda_1(\Omega)$ under geometric assumptions on $\p\Omega$, and to apply them to some shape optimization problems related to $\lambda_1(\Omega)$. 

It turns out to be more convenient when dealing with $\lambda_1(\Omega)$ to rewrite (\ref{BestConst_1}) as a minimization problem in the space $BV(\Omega)$ of functions of bounded variation (see \cite{AFP,EvaGar,Ziemer}) in the following way: 
\begin{equation}\label{BestConst_2}
  \lambda_1(\Omega)=\inf_{u\in BV(\Omega),~u\not\equiv 0~on~\p\Omega} 
 \frac{\displaystyle \int_{\Omega} |\nabla u| + \displaystyle \int_{\Omega}|u|\, dx}{\displaystyle \int_{\partial\Omega} |u|\,dH^{N-1}}. 
\end{equation}
The equivalence between (\ref{BestConst_1}) and (\ref{BestConst_2}) follows from the fact that given $u\in BV(\Omega)$, there exist $u_n\in C^\infty(\Omega)$ such that $u_n=u$ on $\p\Omega$ and the $u_n$'s approximate $u$ in the sense that $u_n\to u$ in $L^1(\Omega)$ and 
$\int_\Omega |\nabla u_n| dx\to \int_\Omega |\nabla u|$ (see \cite{Demengel}, \cite{Giusti}). 

We can also express $\lambda_1(\Omega)$ in a more geometric way as an isoperimetric type problem. We recall that a set $A\subset\bar\Omega$ is said of finite perimeter if its characteristic function $\chi_A$ belongs to $BV(\R^n)$. It then follows from the coarea formula that 
\begin{equation}\label{BestConst_3}
 \begin{split}
  \lambda_1(\Omega)
  & = \inf_{A\subset\bar\Omega,\chi_A\in BV(\R^n)} \frac{\displaystyle |\p A\cap \Omega| + |A|}{\displaystyle |A\cap \p\Omega|}, \\
 \end{split}
\end{equation} 
where $|\p A\cap \Omega|$ and $|A\cap \p\Omega|$ stands for $H^{n-1}(\p A\cap \Omega)$ and $H^{n-1}(A\cap \p\Omega)$ respectively. This infimum is always attained by some set of finite perimeter $A\subset\bar \Omega$ that we call an eigenset. We refer the reader to \cite{LI} for a detailed proof of this result. 

We end this presentation of $\lambda_1(\Omega)$ by recalling its value in the case where $\Omega=B_0(R)$ is a ball or an annulus $\Omega=B_0(R)\backslash \bar B_0(r)$. As remarked in \cite[Remark 1]{AMR}, it follows from \cite{Motron} that 
\begin{equation}\label{Exemples}
 \lambda_1(\Omega)=
 \begin{cases} 
  \frac{\displaystyle |\Omega|}{\displaystyle |\p\Omega|}~\text{ if } 
                                                 \frac{\displaystyle |\Omega|}{\displaystyle |\p\Omega|}\le 1  \\
  1 \hspace{.8cm} \text{ otherwise. }
 \end{cases}
\end{equation}
Moreover, if $|\Omega|/|\p\Omega|\le 1$, then $u=|\p\Omega|^{-1}\chi_\Omega$ is a minimizer, and the only normalized one if $|\Omega|/|\p\Omega|= 1$, whereas if $|\Omega|/|\p\Omega|\ge 1$, there is no extremal for $\lambda_1(\Omega)$. 

\medskip

We first consider the problem of the existence of an extremal for $\lambda_1(\Omega)$. Since the immersion $W^{1,1}(\Omega)\hookrightarrow L^1(\partial\Omega)$ is not compact, the existence of minimizers for $\lambda_1(\Omega)$ does not follows by standard methods. Indeed this problem has already been considered in \cite{AMR} and \cite{Demengel} where it is proved that $\lambda_1(\Omega)$ is attained as soon as 
\begin{equation}\label{Cond}
 \lambda_1(\Omega)<1. 
\end{equation}
We will provide an alternative proof of this result. 
Notice that according to \cite{AMR,Demengel}, the large inequality in (\ref{Cond}) always holds. We refer to \cite{AMR} for the derivation of the Euler equation satisfed by a minimizer. 
According to \cite{Motron}, $\lambda=1$ is the best first constant in the embedding $W^{1,1}(\Omega)\hookrightarrow L^1(\partial\Omega)$ in the sense that for any $\epsilon>0$ there exists $B_\epsilon>0$ such that for any $u\in BV(\Omega) $, 
\begin{equation}\label{BestInequ}
 \int_{\p\Omega} |u|\,dH^{N-1} \le (1+\epsilon)\int_\Omega |\nabla u| + B_\epsilon\int_\Omega |u|\,dx,
\end{equation}
and $1$ is the lowest constant such that such an inequality holds for any $\epsilon>0$ and any $u\in BV(\Omega) $. The inequality (\ref{Cond}) is then the usual condition ensuring that $\lambda_1(\Omega)$ is attained when dealing with critical problem (see e.g. \cite{Aubin1}, \cite{DH}). 

Our first result provides a local geometric condition on $\Omega$ for (\ref{Cond}) to hold. Before stating it, we need a definition. 
We say that a point $x\in\p\Omega$ is a {\it "good point"} if the curvature of $\p\Omega$ at $x$ is big enough, more precisely if the principal curvatures $\lambda_1,\dots,\lambda_{N-1}$ of $\p\Omega$ at $x$ are all positive  and satisfy 
$ \sum_{i=1}^{N-1} \lambda_i > 1$,
and if the graph of $\p\Omega$ around $x$ is close to the parabola $y\to(1/2)\sum \lambda_i y_i^2$ when considered in a local coordinate system such that $x=0$ and the unit outward normal derivative at $0$ of $\p\Omega$ is $(0,\dots,0,1)$ (see (\ref{HypBord}) for a precise statement). 

The result is the following:

\begin{thm}\label{thm1} 
 If there exists a "good point" $x\in\p\Omega$, then (\ref{Cond}) holds.
\end{thm}

\noindent Similarly, we can also prove that (\ref{Cond}) holds when a part of $\p\Omega$ is close to a convex cone of vertex $x\in\p\Omega$ and angle in $(0,\pi/2)$, that is a non-flat cone, since in that case the "curvature" of $\p\Omega$ at $x$ is infinite. 

\medskip

It is well-known that for $p>1$, the trace embedding $W^{1,p}(\Omega)\hookrightarrow L^p(\p\Omega)$ is continuous and compact. In particular the best constant $\lambda_p(\Omega)$ for this embedding, namely
\begin{equation*}
 \lambda_p(\Omega)=\inf_{u\in W^{1,p}(\Omega)\setminus W_0^{1,p}(\Omega)} 
              \frac{\displaystyle \int_{\Omega} |\nabla u|^p + |u|^p\, dx}{\displaystyle \int_{\partial\Omega} |u|^{p}\,dH^{N-1}}, 
\end{equation*}
is attained by some positive $u_p$ normalized by $\int_{\p\Omega} u_p^p\,dH^{N-1} = 1$. 
To show the existence of an extremal for $\lambda_1(\Omega)$, the authors of \cite{AMR} approached $\lambda_1(\Omega)$ by $\lambda_p(\Omega)$. They proved that 
\begin{equation}\label{ConvEigenvalue}
 \lambda_p(\Omega)\to \lambda_1(\Omega)\mbox{ as } p\to 1,
\end{equation}
and also that 

\begin{thm}\label{thm2}
 if $\lambda_1(\Omega)<1$, there exists a nonnegative function $u\in BV(\Omega)$ normalized by $\int_{\p\Omega} |u|\,dH^{N-1} =1$, which attains the
 infimum in the  definition of $\lambda_1(\Omega)$, and such that 
 $$ u_p\to u~\text{ in $L^1(\p\Omega)$ }~\text{ and } ~\int_\Omega |\nabla u_p|^p\,dx \to \int_\Omega |\nabla u| $$
 as $p\to 1$. 
\end{thm} 

\noindent We will give a short proof of this result, different from the one provided in \cite{AMR,Demengel}.  

As an immediate corollary, we have that 

\begin{cor} If $\p\Omega$ has a "good point", then $\lambda_1(\Omega)$ is attained. 
\end{cor}

\medskip

As an application of Theorem \ref{thm1}, we study a shape optimization problem related to $\lambda_1(\Omega)$. Given $\alpha\in (0,|\Omega|)$, where $|\Omega|$ denotes the volume of $\Omega$, and a measurable subset $A\subset\Omega$ of volume $\alpha$, we first consider the minimization problems
\begin{equation*}
 \lambda_{1,A} = \inf_{\begin{cases} u\in BV(\Omega),~u\not\equiv 0~on~\p\Omega\\ u=0~in~A \end{cases}}
 \frac{\displaystyle \int_\Omega |\nabla u| + |u|\,dx}{\displaystyle \int_{\p\Omega} |u|\,dH^{N-1}},
\end{equation*}
and
\begin{equation*}
 \lambda_{p,A} = \inf_{\begin{cases} u\in W^{1,p}(\Omega)\setminus W^{1,p}_0(\Omega)\\ u=0~in~A \end{cases}}
 \frac{\displaystyle \int_\Omega |\nabla u|^p + |u|^p\,dx}{\displaystyle \int_{\p\Omega} |u|^p\,dH^{N-1}}.
\end{equation*}

\noindent It is easily seen that $\lambda_{p,A}$, $p>1$, is attained. Concerning $\lambda_{1,A}$, we have, in the same spirit as what we had for $\lambda_1(\Omega)$, that 

\begin{thm}\label{thmHole}
 If 
 \begin{equation*}
  \lambda_{1,A}<1, 
 \end{equation*} 
 there exists an extremal for $\lambda_{1,A}$. Moreover this inequality holds as soon as there exists a good point $x\in\p\Omega$ such that $A\cap
 B_x(r)=\emptyset$ for some $r>0$.
\end{thm}

Remark that $\lambda_{p,A}$, $p\ge 1$, does not change if we modify $A$ on a set of Lebesgue measure zero. To give a meaning to $\lambda_{p,A}$, $p>1$, when $|A|=0$, the authors of \cite{BRW} modified $\lambda_{p,A}$ by minimizing over $\overline{C^\infty_c(\bar\Omega\backslash A)}$. In the case $p=1$, we introduce in a similar way the set $BV_A(\Omega)$ of the functions $u\in BV(\Omega)$ that can be approximated by a sequence $u_\epsilon\in C^\infty_c(\bar\Omega\backslash A)$ in the sense that $u_\epsilon\to u$ in $L^1(\Omega)$ and 
$\int_\Omega |\nabla u_\epsilon|\to \int_\Omega |\nabla u|$. We can then prove as in \cite{EvaGar} that $BV_A(\Omega)=BV(\Omega)$ if and only if 
$\text{cap}_1(A)=0$, where $ \text{cap}_1(A)$ denotes the $1$-capacity of $A$ defined by 
$$ \text{cap}_1(A)= \inf~\left\{\int_{\R^n} |\nabla u|,~~u\in BV(\R^n),~A\subset\text{int}\{u\ge 1\}\right\}.$$
In the case where $A$ is compact, the coarea formula implies that $\text{cap}_1(A)=\inf~|\p\omega|$ where the infimum is taken over all the smooth open subsets $\omega\subset\R^n$ containing $A$ (see \cite{M}). We consider the minimization problem 
$$ \lambda_{1,A}'= \inf_{u\in BV_A(\Omega)}
 \frac{\displaystyle \int_\Omega |\nabla u| + |u|\,dx}{\displaystyle \int_{\p\Omega} |u|\,dH^{N-1}}.$$
Then $\lambda_{1,A}\le \lambda_{1,A}'$ with equality when $ \text{cap}_1(A)=0$. If $ \text{cap}_1(A)>0$, both cases $\lambda_{1,A}=\lambda_{1,A}'$ and 
$\lambda_{1,A}< \lambda_{1,A}'$ can occur. For example if a part of the boundary of $\Omega\subset\R^2$ has curvature big enough (e.g. like a smooth version of the set $Q_{\delta,\eta}$ defined below next to theorem \ref{thm4}), then $\lambda_1(\Omega)$ will be attained by some $\chi_C$ where $C\subsetneq \Omega$. Then if we put a small curve $A$ in the interior of $\Omega\backslash C$, $\chi_C\in BV_A(\Omega)$ and thus $\lambda_\emptyset=\lambda_{1,A}=\lambda_{1,A}'$. On the contrary, if $\Omega\subset\R^2$ is a ball such that $|\p\Omega|=|\Omega|$, then we know that $\lambda_1(\Omega)$ is attained only by the $\mu\chi_\Omega$, $\mu\in\R$. Then if $A$ small segment inside $\Omega$,  $\lambda_{1,A}< \lambda_{1,A}'$. 

\medskip

We now want to minimize $\lambda_{p,A}$, $p\ge 1$, when $A$ runs over all the measurable subsets of $\Omega$ of volume $\alpha$ i.e. we look at the following shape optimization problem:
$$ \lambda_p(\alpha)= \inf_{A\subset\Omega,~|A|=\alpha} \lambda_{p,A} $$
for $p\ge 1$ and $\alpha\in (0,|\Omega)$. 

The optimization problem $\lambda_p(\alpha)$, $p>1$, has been considered recently. 
Existence of an optimal set  has been established in \cite{BRW}, and its regularity investigated in \cite{BRW2} for $p=2$. The optimization problem $\lambda_p(\alpha)$ with a critical exponent has been considered in \cite{BS}. 
Such problems of optimal design appear in several branches of applied mathematics, specially in the case $p=2$. For example in problems of minimization of the energy stored in the design under a prescribed loading. We refer to \cite{CC} for more details.
 
We prove the following relation between $\lambda_p(\alpha)$ and $\lambda_1(\alpha)$: 
 
\begin{thm}\label{thm3} 
 We have  
 \begin{equation}\label{prop2_1}
  \limsup_{p\to 1} \lambda_p(\alpha) \le \lambda_1(\alpha). 
 \end{equation}
 Moreover, if there exists a {\it good point} $x\in\p\Omega$, then 
 \begin{equation}\label{prop2_2}
   \lim_{p\to 1} \lambda_p(\alpha) = \lambda_1(\alpha). 
 \end{equation}
\end{thm}

\noindent The proof of this theorem gives the existence of an extremal $u\in BV(\Omega)$ for $\lambda_1(\alpha)$ but, since we can only prove that $|\{u=0\}|\ge\alpha$ and not $|\{u=0\}|=\alpha$, we cannot assert the existence of an optimal hole $A$ such that $\lambda_1(\alpha)=\lambda_{1,A}$. However if we consider the following modified optimal design problem 
\begin{equation}\label{Modified}
 \tilde \lambda_1(\alpha)= \inf_{\begin{cases} u\in BV(\Omega),~u\not\equiv 0~on~\p\Omega\\ |\{u=0\}|=\alpha \end{cases}}
 \frac{\displaystyle \int_\Omega |\nabla u| + |u|\,dx}{\displaystyle \int_{\p\Omega} |u|\,dH^{N-1}},
\end{equation}
we can prove that 

\begin{thm}\label{thm3Modified}
 if there exists a {\it good point} $x\in\p\Omega$, then $\tilde \lambda_1(\alpha)$ is attained by some $u$. In particular $\{u=0\}$ is an optimal hole for $\tilde \lambda_1(\alpha)$. 
\end{thm}

\noindent It follows from \cite{BRW} that $\lambda_p(\alpha)=\tilde \lambda_p(\alpha)$, $p>1$, where $\tilde \lambda_p(\alpha)$ is defined by
$$ \tilde \lambda_p(\alpha)= \inf_{\begin{cases} u\in W^{1,p}(\Omega)\setminus W^{1,p}_0(\Omega)\\ |\{u=0\}|=\alpha \end{cases}}
 \frac{\displaystyle \int_\Omega |\nabla u|^p + |u|^p\,dx}{\displaystyle \int_{\p\Omega} |u|^p\,dH^{N-1}},$$
but for the same reason as before, we cannot establish the convergence of $ \tilde \lambda_p(\alpha)$ to $\tilde \lambda_1(\alpha)$ as $p\to 1$ . 

\medskip

Our last result concerning $\lambda_1$ is the computation of the first variation, the so-called shape derivative, of the functional 
$\Omega\to \lambda_1(\Omega)$. 
Let $R:\R^n\to\R^n$ be a $C^1$ vector-field, and $\Omega_\delta = T_\delta(\Omega)$, where $T_\delta$ is the $C^1$-diffeomorphism defined for $\delta$ small by 
 $$T_\delta(x) = x + \delta R(x). $$
We will prove that the map $\delta\to\lambda_1(\Omega_\delta)$ is continuous at $\delta=0$, and also differentiable at $\delta=0$ under an aditional  uniqueness assumption holding for example when $\Omega$ is a ball. 

Remark that if we allow perturbations of the domains that are less regular, we may not have continuity of $\lambda_1(\Omega_\delta)$ as the following example shows. Let $Q=[0,1]^N$ be the unit cube of $\R^N$, and let $Q_{\delta,\eta}=Q\cup A_{\delta,\eta}$ with 
$$ A_{\delta,\eta} = [1,1+\eta]\times [0,\delta]\times [0,1]^{N-2},~\delta,\eta>0.$$
Then taking $\chi_A$ as a test-function to estimate $\lambda_1(Q_\delta)$, we get 
$$ \lambda_1(Q_\delta) \le \dfrac{\delta+\eta\delta}{C\eta}\to 0 $$
as $\delta\to 0$ if $\eta>>\delta$. This shows that, even if $|Q_\delta\Delta Q|\to 0$ or $Q_\delta\to Q$ in Hausdorff distance, we don't have continuity of $\lambda_1(Q_\delta)$. Indeed $\lambda_1(Q_\delta)\to 0\neq \lambda_1(Q)$. 

Shape analysis is the subject of an intense research activity. We refer for example to \cite{Henrot} for an introduction to this field.
To the best of the author's knowledge, the shape analysis of a problem involving the $L^1$-norm of the gradient has only been considered up to know in \cite{HS,Saintier} where the authors deal with the best constant for the embedding of $W^{1,1}(\Omega)$ into $L^1(\Omega)$. 

Our result is the following:   

\begin{thm}\label{thm4}
We have 
$$ \lambda_1(\Omega_\delta)\to \lambda_1(\Omega) $$
as $\delta\to 0$.  
Moreover, if we assume that $\lambda_1(\Omega)<1$ and that there exists a unique nonnegative extremal $u\in BV(\Omega)$ for $\lambda_1(\Omega)$ normalized by $\int_{\p\Omega}  u\,dH^{N-1}=1$, then $u=|A\cap \p\Omega|^{-1}\chi_A$ for some set of finite perimeter $A\subset\bar\Omega$, and the map $\delta\to\lambda_1(\Omega_\delta)$ is differentiable at $\delta=0$ with 
\begin{equation}\label{ShapeDer}
\begin{split}
 & \frac{d}{d\delta}\lambda_1(\Omega_\delta)_{|\delta=0} = \\
 & \int_{\bar\Omega}  \left\{ f(\nu)\chi_{\p^*A\cap\Omega} - \lambda_1(\Omega)f(\vec{n})\chi_{A\cap\p\Omega} - (R,\nu)\chi_{\p^*A}
   \right\}\,\frac{dH^{N-1}}{|A\cap \p\Omega|},
 \end{split}
\end{equation} 
where $f(X)=div~R-(X;DR.X)$, $X\in\R^n$, 
$\nu$ is the Radon-Nikodym derivative of $|\nabla u|$ with respect to $\nabla u$, $\vec{n}$ is the unit outward normal to $\p\Omega$, and $\p^*A$ is the reduced boundary of $A$ (see \cite{AFP,EvaGar,Ziemer}). 
\end{thm}

\noindent As said previously (see the comments next to (\ref{Exemples})), the uniqueness property used in this theorem holds in particular when $\Omega$ is a ball such that $|\p\Omega|=|\Omega|$. Its unique eigenset is then $\bar\Omega$ and $\lambda_1(\Omega)=1$, so that we can rewrite (\ref{ShapeDer}) as 
$$ \frac{d}{d\delta}\lambda_1(\Omega_\delta)_{|\delta=0} 
  = \int_{\p\Omega} \left\{(R,\vec{n})-(div~R-(\vec{n};DR.\vec{n}))\right\}\,\frac{dH^{N-1}}{|\p\Omega|}. $$
Denoting by $div_g$ the divergence operator of the manifold $(\p\Omega,g)$, where $g$ is the metric induced by the Euclidean metric on $\p\Omega$, by $H$ the mean curvature of $\p\Omega$, and by $R_{\p\Omega}$ the tangential part of $R$, we have (see \cite{Henrot}): 
$$ div\,R-(\vec{n};DR.\vec{n}) = div_g\,R_{\p\Omega} + H(R,\vec{n}). $$
Since $\Omega$ is of radius $n$, $H=1/n$, and the previous formula becomes
$$ \frac{d}{d\delta}\lambda_1(\Omega_\delta)_{|\delta=0} 
  = \int_{\p\Omega} (1-H)(R,\vec{n})\,\frac{dH^{N-1}}{|\p\Omega|} 
  = -\frac{n-1}{n}\int_{\p\Omega} (R,\vec{n})\,\frac{dH^{N-1}}{|\p\Omega|}. $$
In particular, if we consider measure-preserving deformation, i.e. vector-fields $R$ such that $div\,R=0$, we get 
$$ \frac{d}{d\delta}\lambda_1(\Omega_\delta)_{|\delta=0} = 0, $$
so that a ball of $\R^n$ of radius $n$ is critical for such deformations. 

\bigskip

The paper is organized as follow. We prove theorem \ref{thm1} - \ref{thm3Modified} in the following section and theorems \ref{thm4} in the last one.

\section{Proof of theorems \ref{thm1} - \ref{thm3}}

\subsection{Proof of theorem \ref{thm1}}
Let $x_0\in\p\Omega$ be a ``good point''. By taking an appropriate coordinate system, we can assume that $x_0=0$ and that there exist $r>0$ such that
\begin{align*}
B_r\cap\Omega = & \{(y,t)\in B_r,\ t>\rho(y)\} \\
B_r\cap\p\Omega =&  \{(y,t)\in B_r,\ t=\rho(y)\}
\end{align*}
where $y=(y_1,\dots,y_{N-1})\in\R^{N-1}$, $B_r$ is the Euclidean ball centered at the origin and of radius $r$, and
$$ \rho(y)= \frac{1}{2}|y|_\lambda^2(1+O(|y|^\alpha)) $$
for some $\alpha>0$, with 
$$ |y|_\lambda^2 = \sum_{i=1}^{N-1} \lambda_i y_i^2, $$
where the $\lambda_i$'s are the principal curvatures of $\p\Omega$ at $0$. 
We assume that $\alpha$ is such that as $\epsilon\to 0$, 
\begin{equation}\label{HypBord}
 |\{y\in\R^{N-1},~~\rho(y)\le \epsilon^2/2\}\Delta \{y\in\R^{N-1},~~|y|_\lambda\le \epsilon\}| = o(\epsilon^{N+1}),
\end{equation}
where $A\Delta B = (A\setminus B)\cup (B\setminus A)$ denotes the symetric difference of the sets $A,B\subset\R^{N-1}$ and $|A|$ the volume of $A$.
A sufficient condition for (\ref{HypBord}) to hold is $\alpha>2$. 

We consider the test-functions 
$$ u_\epsilon(y,t)=\chi_{\Omega\cap \{0\le t\le\epsilon^2/2\}}(y,t). $$
Assume for the moment that the following asymptotic developments hold:  
\begin{equation}\label{NormeGradient2} 
   \int_\Omega |\nabla u_{\epsilon}|  =  ~b_{N-1}^\lambda\epsilon^{N-1} + o(\epsilon^{N+1}),  
\end{equation}
\begin{equation}\label{NormeL^p2} 
   \int_\Omega |u_{\epsilon}|\,dydt  =  \frac{\omega_{N-2}^\xi}{2(N+1)(N-1)\sqrt{\prod \lambda_i}}~\epsilon^{N+1} + o(\epsilon^{N+1}), 
\end{equation}
and
\begin{equation}\label{NormeBord2} 
   \int_{\p\Omega} |u_{\epsilon}|\,dH^{N-1}  =  ~\epsilon^{N-1}b_{N-1}^\lambda 
   + \frac{\omega_{N-2}^\xi\sum \lambda_i}{2(N-1)(N+1)\sqrt{\prod\lambda_i}}~ \epsilon^{N+1} + o(\epsilon^{N+1}),
\end{equation}
where $b_{N-1}^\lambda = |\{y\in\R^{N-1},~~|y|_\lambda\le 1\}|$ and $\omega_{N-2}^\xi = |\{y\in\R^{N-1},~~\sum y_i^2=1\}|$.
It then follows that 
\begin{equation*}
\begin{split}
 \lambda_1 & \le  \frac{\displaystyle \int_{\Omega} |\nabla u_\epsilon| + \displaystyle \int_{\Omega}|u_\epsilon|\, dx}
         {\displaystyle \int_{\partial\Omega} |u_\epsilon|\,dH^{N-1}} \\
 & = 1 + \frac{\omega_{N-2}^\xi}{2(N-1)(N+1)b_{N-1}^\lambda\sqrt{\prod\lambda_i}} \left\{1-\sum\lambda_i\right\}\epsilon^2 + o(\epsilon^2), 
\end{split}
\end{equation*}
from which we deduce Theorem \ref{thm1}. 

We now prove (\ref{NormeGradient2}), (\ref{NormeL^p2}) and (\ref{NormeBord2}). In view of (\ref{HypBord}), 
\begin{equation*}
\begin{split}
 \int_\Omega |\nabla u_\epsilon|
  = & ~|\{\rho(y)\le\epsilon^2/2\}| = ~|\{|y|_\lambda\le\epsilon\}| + o(\epsilon^{N+1}) \\
  = & ~\epsilon^{N-1}b_{N-1}^\lambda + o(\epsilon^{N+1}) 
\end{split}
\end{equation*}
which proves (\ref{NormeGradient2}). We now prove (\ref{NormeL^p2}). We first note that 
\begin{equation*}
\begin{split}
 \int_\Omega |u_\epsilon|\,dydt 
  = & \int_{\{\rho(y)\le\epsilon^2/2\}} \left(\int_{\rho(y)}^{\epsilon^2/2} dt\right) dy \\
  = & ~\frac{\epsilon^2}{2}|\{|y|_\lambda\le\epsilon\}| -\int_{\{|y|_\lambda\le\epsilon\}} \frac{1}{2}|y|_\lambda^2(1+O(|y|^\alpha))\,dy 
       + o(\epsilon^{N+1}) \\
  = &  ~\frac{b_{N-1}^\lambda}{2}\epsilon^{N+1} - \frac{\epsilon^{N+1}}{2}\int_{ \{|y|_\lambda\le 1\} } |y|_\lambda^2\,dy + o(\epsilon^{N+1}). 
\end{split}
\end{equation*}
Denoting by $b_{N-1}^\xi$ (resp. $\omega_{N-2}^\xi$) the volume of the unit ball (resp. the unit sphere) of $\R^{N-1}$ for the usual Euclidean metric $\xi$, we have 
$$ b_{N-1}^\lambda = \frac{b_{N-1}^\xi}{\sqrt{\prod \lambda_i}} = \frac{\omega_{N-2}^\xi}{(N-1)\sqrt{\prod \lambda_i}}, $$
and, by the coarea formula, 
\begin{equation*}
\begin{split}
 \int_{ \{|y|_\lambda\le 1\} } |y|_\lambda^2\,dy 
 & = \frac{1}{\sqrt{\prod \lambda_i}} \int_0^1 \left(\int_{\{|y|_\xi=t\}} |y|^2_\xi\,dH^{N-2}\right)\,dt \\
 &  = \frac{\omega_{N-2}^\xi}{(N+1)\sqrt{\prod \lambda_i}}.
\end{split}
\end{equation*}
Hence 
$$ \int_\Omega |u_\epsilon|\,dydt = \frac{\omega_{N-2}^\xi}{2(N+1)(N-1)\sqrt{\prod \lambda_i}}\epsilon^{N+1} + o(\epsilon^{N+1}) $$
which is (\ref{NormeL^p2}). Eventually, to prove (\ref{NormeBord2}), we write that 
\begin{equation*}
\begin{split}
 \int_{\p\Omega} |u_\epsilon|\,dH^{N-1} 
 = & \int_{\{\rho(y)\le\epsilon^2/2\}} \sqrt{1+|\nabla \rho|^2}\, dy \\
 = & \int_{\{|y|_\lambda\le\epsilon\}} \sqrt{1+|\nabla \rho|^2}\, dy + o(\epsilon^{N+1}) \\
 = & \int_{\{|y|_\lambda\le\epsilon\}} (1+\frac{1}{2}\sum \lambda_i^2 y_i^2 + o(|y|_\lambda^2))\, dy + o(\epsilon^{N+1})  \\
 = & ~\epsilon^{N-1}b_{N-1}^\lambda + \frac{\epsilon^{N+1}}{2} \int_{\{|y|_\lambda\le 1\}} \sum \lambda_i^2 y_i^2 dy + o(\epsilon^{N+1}) \\
\end{split}
\end{equation*}
with, using the symetry of the sphere and then the coarea formula,  
\begin{equation*}
\begin{split}
 \int_{\{|y|_\lambda\le 1\}} \sum \lambda_i^2 y_i^2 dy 
   & = \frac{\sum \lambda_i}{\sqrt{\prod\lambda_i}}  \int_{\{|y|_\xi\le 1\}} y_i^2\, dy \\
   & = \frac{\sum \lambda_i}{(N-1)\sqrt{\prod\lambda_i}}  \int_{\{|y|_\xi\le 1\}} |y|_\xi^2\, dy  \\
   & = \frac{\omega_{N-2}^\xi\sum \lambda_i}{(N-1)(N+1)\sqrt{\prod\lambda_i}}. 
\end{split}
\end{equation*}
Hence 
\begin{equation*}
\begin{split}
 \int_{\p\Omega} |u_\epsilon|\,dH^{N-1} 
  = & ~\epsilon^{N-1}b_{N-1}^\lambda + \frac{\omega_{N-2}^\xi\sum \lambda_i}{2(N-1)(N+1)\sqrt{\prod\lambda_i}} \epsilon^{N+1} + o(\epsilon^{N+1}) \\
\end{split}
\end{equation*}
which is (\ref{NormeBord2}). 
 
\medskip

We now assume that, at a point $x\in\p\Omega$, $\Omega$ is close to the cone $C_\omega=\{\lambda\omega,~\lambda\ge 0\}$, where $\omega$ is a subset of the unit sphere of $\R^N$, in the sense that 
\begin{eqnarray*}
 |\epsilon^{-1}(\Omega-x)\cap B_0(1)| & \sim & |C_\omega\cap B_0(1)|,\\
 |\epsilon^{-1}\p(\Omega-x)\cap B_0(1)| & \sim & |\p C_\omega\cap B_0(1)|, \\
 |\epsilon^{-1}(\Omega-x)\cap \p B_0(1)| & \sim & |C_\omega\cap \p B_0(1)| 
\end{eqnarray*}
as $\epsilon\to 0$. Using $u_\epsilon=\chi_{\Omega\cap B_x(\epsilon)}$ as a test-function, we have 
\begin{eqnarray*}
 \int_\Omega u_\epsilon\, dx & = & |\Omega\cap B_x(\epsilon)|\sim \epsilon^N|C_\omega\cap B_0(1)|,\\
 \int_{\p\Omega} u_\epsilon\, d\sigma & = & |\p\Omega\cap B_x(\epsilon)| \sim \epsilon^{N-1} |\p C_\omega\cap B_0(1)|,\\
 \int_{\p\Omega} |\nabla u_\epsilon| & = & |\Omega\cap \p B_x(\epsilon)|\sim \epsilon^{N-1}|C_\omega\cap \p B_0(1)|= \epsilon^{n-1}|\omega|,
\end{eqnarray*}
with 
$$ |\p C_\omega\cap B_0(1)| = \int_0^1 |\p (r\omega)|\,dr = \frac{|\p\omega|}{N-1}, $$
and thus 
$$ \lambda_1 \le \frac{|\omega|}{|\p C_\omega\cap B_0(1)|}+O(\epsilon) = \frac{(N-1)|\omega|}{|\p\omega|}+O(\epsilon).$$
Hence if $(N-1)|\omega|< |\p\omega|$, we get (\ref{Cond}).
In the particular case where $\omega$ is a spherical cap, i.e. the intersection of $\p B_0(1)$ with an half-space $H^+$ defined by an affine hyperplane $H$, in such a way that $C_\omega$ is convex of angle $\alpha\in (0,\pi/2]$, we can get in a similar way that 
\begin{eqnarray*}
 \lambda_1 & \lesssim & \frac{(N-1)|H\cap B_0(1)|}{|H\cap \p B_0(1)|}=\frac{(N-1)\sin^{N-1}(\alpha)b_{N-1}^\xi}{\sin^{N-2}(\alpha)\omega_{N-2}^\xi} \\
           & = & \sin(\alpha).
\end{eqnarray*}
Hence if $\epsilon^{-1}(\Omega-x)$ is asymptotically close to the cone $C_\omega$ with angle $\alpha\in (0,\pi/2)$, (\ref{Cond}) holds.

\subsection{Proof of theorem \ref{thm2}}

We adapt to our case the argument of \cite{Demengel2}. 
In view of (\ref{ConvEigenvalue}), the sequence $(\lambda_p)_{p>1}$ is bounded, from which it follows that the sequence $(\|u_p\|_{W^{1,p}})$ is bounded, and eventually that the sequence $(u_p)$ is bounded in $BV(\Omega)$. In particular, there exists $u\in BV(\Omega)$ such that, up to a subsequence, $u_p\to u$ strongly in $L^q(\Omega)$ for all $q<N/(N-1)$ and a.e.. In particular, $u\ge 0$ a.e..
According to \cite{Lions} (see also \cite{Demengel}) and in view of (\ref{BestInequ}), there exist a nonempty set $I\subset\mathbb{N}$, a sequence of points $(x_i)_{i\in I}\subset\p\Omega$ and sequences of positive reals $(\mu_i)_{i\in I}$, $(\nu_i)_{i\in I}$, and two measures $\mu$ and $\nu$, with $supp~\nu\subset\p\Omega$, such that 
\begin{equation}\label{CCP}
\begin{cases}
\begin{split}
 |\nabla u_p|^p\,dx & \rightharpoonup \mu \ge |\nabla u| + \sum_{i\in I} \nu_i \delta_{x_i}, \\ 
 |u_p|^p\,dH^{N-1} & \rightharpoonup \nu= |u|\,dH^{N-1} + \sum_{i\in I} \nu_i \delta_{x_i}. 
\end{split}
\end{cases}
\end{equation}
Let $\sigma_p=|\nabla u_p|^{p-2}\nabla u_p$. Given $q\in [1,+\infty)$, it is easily seen, using Hölder' inequality, that $(\sigma_p)$ is bounded in $L^q(\Omega)$ for $p$ small enough. Hence there exists $\sigma\in \cap_{q\ge 1} L^q(\Omega)$ such that $ \sigma_p\to \sigma$ weakly in $L^q(\Omega)$ for every $q>1$. Notice that $\sigma\in L^\infty(\Omega)$ with $\|\sigma\|_\infty\le 1$. Indeed for any $\psi\in C_c^\infty(\Omega,\R^n)$, we have 
$$ \left|\int_\Omega \sigma\psi\,dx\right| = \lim_{p\to 1} \left|\int_\Omega \sigma_p\psi\,dx\right| 
 \le \lim_{p\to 1} \|\nabla u_p\|_p^{p-1}\|\psi\|_p = \int_\Omega |\psi|\,dx. $$
Passing to the limit in the Euler equation for $u_p$, namely
\begin{equation}\label{EulerEqu}
 \int_\Omega \sigma_p\nabla\psi\,dx + \int_\Omega u_p^{p-1}\psi\,dx = \lambda_p(\Omega) \int_{\p\Omega} u_p^{p-1}\psi\,dH^{N-1},~\hspace{.1cm}
\forall~\psi\in W^{1,p}(\bar\Omega), 
\end{equation}
we get that, in view of (\ref{ConvEigenvalue}), that  
\begin{equation}\label{Equ_sigma} 
\begin{cases}
 -div~\sigma+1=0~\hspace{1cm} \text{ in }\Omega \\
 \hspace{1.3cm} \sigma.\vec{n}=\lambda_1(\Omega)~\hspace{0.25cm} \text{ on }~\p\Omega,
\end{cases}
\end{equation}
where $\vec{n}$ is the unit outward normal to $\p\Omega$. 
Let $\phi\in C^\infty(\bar\Omega)$. Passing to the limit in (\ref{EulerEqu}) with $\psi=u_p\phi$, using (\ref{ConvEigenvalue}), we obtain 
\begin{equation}\label{Equ1}
 \int_\Omega \phi\,d\mu + \int_\Omega u\sigma\nabla\phi\,dx + \int_\Omega u\phi\,dx = \lambda_1(\Omega)\int_{\p\Omega} \phi\,d\nu.  
\end{equation}
According to the definition of the measure $\sigma\nabla u$, defined weakly by integration by part (see \cite{Demengel2}), and in view of (\ref{Equ_sigma}), we have 
\begin{equation}\label{Equ2}
\begin{split}
 \int_\Omega u\sigma\nabla\phi\,dx 
 & = \int_\Omega div~(\phi u\sigma)\,dx - \int_\Omega \phi u (div~\sigma)\,dx - \int_\Omega \phi (\sigma\nabla u) \\
 & = \lambda_1(\Omega) \int_{\p\Omega} \phi u\,dH^{N-1} -\int_\Omega \phi u\,dx - \int_\Omega \phi (\sigma\nabla u). 
\end{split}
\end{equation} 
Plugging this in (\ref{Equ1}) and using the definition of $\mu$ and $\nu$, we eventually get 
$$ \int_\Omega \phi (|\nabla u|-\sigma\nabla u) \le (\lambda_1-1) \int_\Omega \phi (\sum_{i\in I} \nu_i \delta_{x_i}). $$
Since $|\sigma\nabla u|\le \|\sigma\|_\infty|\nabla u|\le |\nabla u|$ and $\lambda_1<1$ by assumption, we deduce that $\nu_i=0$ for all $i\in I$. 
In particular $\int_{\p\Omega} u\,dH^{N-1}=1$. Moreover, inserting (\ref{Equ2}) into (\ref{Equ1}), we see that $\mu = \sigma\nabla u \le |\nabla u|$. Hence $\mu=|\nabla u|$.

\subsection{Proof of theorem \ref{thmHole}} The proof of the first part is analogous to the proof of theorem \ref{thm2}. Concerning the second part, just remark that since the principal curvatures at a the {\it good point} $x\in\p\Omega$ are positive, we have 
$supp~~u_\epsilon\subset B_x(r)$ for $\epsilon$ small, where $u_\epsilon$ is the sequence of test-functions considered in the proof of theorem \ref{thm1}. Hence the $u_\epsilon$'s are also admissible test-functions for $\lambda_{1,A}$.

\subsection{Proof of theorem \ref{thm3}}  
 We first prove (\ref{prop2_1}). Given $\epsilon>0$, let $D\subset \Omega$ measurable, $|D|=\alpha$, be such that 
 $$ \lambda_1(D) \le \lambda_1(\alpha) + \epsilon.$$
 The same arguments used to prove (\ref{ConvEigenvalue}) shows that $\lambda_p(D)\to \lambda_1(D)$ as $p\to 1$ (see \cite{AMR}). Hence 
 $$ \limsup_{p\to 1} \lambda_p(\alpha) \le \lim_{p\to 1} \lambda_p(D) = \lambda_1(D) \le \lambda_1(\alpha) + \epsilon.$$
 Since $\epsilon$ is arbitrary, we deduce (\ref{prop2_1}).
 
 Concerning (\ref{prop2_2}), we first note that 
 \begin{equation*}
 \begin{split}
  \lambda_p(\alpha) = \inf_{u\in W^{1,p}(\Omega),|\{u=0\}|\ge\alpha}  
                   \frac{\displaystyle\int_\Omega |\nabla u|^p + |u|^p\, dx }{\displaystyle\int_{\p\Omega} |u|^p\, dH^{N-1} }. 
 \end{split}
 \end{equation*}
 and, in the same way, 
 $$ \lambda_1(\alpha) = \inf_{u\in BV(\Omega),|\{u=0\}|\ge\alpha}  
                   \frac{\displaystyle\int_\Omega |\nabla u| + \int_\Omega |u|\, dx }{\displaystyle\int_{\p\Omega} |u|\, dH^{N-1} }. $$
 For $p>1$, it is known (see \cite{BRW}) that the last infimum is attained by some non-negative $u_p$ normalized by 
 $\int_{\p\Omega} |u_p|^p dH^{N-1}=1$, and satisfying $|\{u_p=0\}|=\alpha$. 
 Independently, since there exists a {\it good point} $x\in\p\Omega$, we have 
\begin{equation}\label{remark}
 \lambda_1(\alpha) < 1.  
\end{equation} 
Indeed, let $D\subset\Omega$ measurable of volume $\alpha$ and consider $D':=(D\backslash B_x(r))\cup \bar D$ for a small $r>0$ and 
$\bar D\subset\Omega$ being such that $|D'|=\alpha$ and $\bar D\subset\Omega\backslash B_x(r)$. Then $D'\cap B_x(r)=\emptyset$, and  thus, according to theorem \ref{thm1}, 
$$ \lambda_1(\alpha)\le  \lambda_1(D') < 1, $$ 
as we wanted to prove. Now, as in the proof of theorem \ref{thm1} and in view of (\ref{remark}), we have that, along a subsequence, 
 \begin{equation*}
 \begin{cases}
  u_p^p  \to u~\text{in } L^1(\Omega)~\text{ and } a.e. \\
  \int_\Omega |\nabla u_p|^p \,dx  \to  \int_\Omega |\nabla u| \\
  \int_{\p\Omega} u\,dH^{N-1}   = \lim_{p\to 1} \int_{\p\Omega} u_p^p\,dH^{N-1} = 1
 \end{cases}
 \end{equation*}
 as $p\to 1$, for some non-negative $u\in BV(\Omega)$. In particular $|\{u=0\}|\ge\alpha$. Hence 
 \begin{equation*}
 \begin{split}
  \lambda_p(\alpha) & = \int_\Omega |\nabla u_p|^p \,dx + \int_\Omega |u_p|^p \,dx = \int_\Omega |\nabla u| + \int_\Omega |u| \,dx + o(1) \\
              & \ge \lambda_1(\alpha).
 \end{split} 
 \end{equation*}
This proves (\ref{prop2_2}). 

\subsection{Proof of theorem \ref{thm3Modified}} A straightforward modification of the proof of (\ref{BestConst_3}) allows us to rewrite (\ref{Modified}) as  
\begin{equation}\label{tilde}
 \tilde \lambda_1(\alpha)=\inf_{\begin{cases} C\subset\bar\Omega,\chi_C\in BV(\R^n) \\ |\Omega\backslash C|=\alpha \end{cases}} 
  \frac{|\p C\cap\Omega|+|C|}{|C\cap\p\Omega}.
\end{equation}
Let $(C_n)$ be a minimizing sequence for this problem. As in the proof of Theorem \ref{thm3}, the existence of a {\it good point} $x\in\p\Omega$ implies that 
\begin{equation}\label{Cond2}
\tilde \lambda_1(\alpha)<1. 
\end{equation}
In particular, for $n$ large enough,
$$ |\p C_n\cap\Omega|+|C_n|\le 2|C_n\cap\p\Omega|\le 2|\p\Omega|,$$
from which we deduce that $(\chi_{C_n})$ is bounded in $BV(\Omega)$. Hence there exists a set of finite perimeter $C$ such that $\chi_{C_n}\to \chi_C$ in $L^1(\Omega)$ and a.e.. In particular $|\Omega\backslash C|=\alpha$. Moreover, as in the proof of theorem \ref{thm4} below, we can deduce from (\ref{Cond2}) that $\int_\Omega |\nabla \chi_{C_n}|\to \int_\Omega |\nabla \chi_{C_n}|$, 
i.e. $|\p C_n\cap\Omega|\to |\p C\cap\Omega|$, and $\int_{\p\Omega|} \chi_{C_n}\,dH^{N-1}\to \int_{\p\Omega|} \chi_{C}\,dH^{N-1}$, i.e. $|C_n\cap\p\Omega|\to |C_n\cap\p\Omega|$. Hence $C$ attains the infimum in (\ref{tilde}), which proves Theorem \ref{thm3Modified}.



\section{Proof of theorems \ref{thm4}}

To simplify the notation, we let $\lambda= \lambda_1(\Omega)$ and $\lambda_\delta=\lambda_1(\Omega_\delta)$. 

According to the change of variable formula for functions of bounded variations \cite{Giusti}, and the change of variable formula for the boundary integral \cite{Henrot}, we have that  
$$ \lambda_\delta = \inf_{u\in BV(\Omega),~u\not\equiv 0~on~\p\Omega} Q_\delta(u) $$
with 
$$ Q_\delta(u)= \frac{ \displaystyle \int_\Omega |(DT_\delta)^{-1}\nu||det~DT_\delta|\,|\nabla u| 
                      + \displaystyle\int_\Omega |u| |det~DT_\delta|\,\,dx }
                     { \displaystyle \int_{\p\Omega} |u| |^t (DT_\delta)^{-1}\vec{n}||det~DT_\delta|\,dH^{N-1} }, $$
where $\nu$ is the Radon-Nikodym derivative of $|\nabla u|$ with respect to $\nabla u$, and $\vec{n}$ is the unit outward normal to $\Omega$.  
We also let $Q=Q_0$, namely 
$$ Q(u) = \frac{ \displaystyle \int_\Omega |\nabla u| + \displaystyle \int_\Omega |u| \,dx }
               {\displaystyle \int_{\p\Omega} |u|\,dH^{N-1} }, $$                     
so that 
$$ \lambda_\delta = \inf_{u\in BV(\Omega),~u\not\equiv 0~on~\p\Omega} Q(u). $$
We first prove that for any $u\in BV(\Omega)$,  
$$ Q_\delta(u)=(1+O(\delta))Q(u) $$
where the $O(\delta)$ is uniform in $u$. The continuity of $\delta\to \lambda_\delta$ at $\delta=0$ then easily follows. 
Let $u\in BV(\Omega$. Since $|\nu|=1$ $|\nabla u|$-a.e., we can assume that $|\nu|=1$ everywhere. Then 
\begin{equation}\label{Asympt1}
 |(DT_\delta)^{-1}\nu| = 1-(\nu,DR.\nu)\delta + o(\delta),
\end{equation}
and in the same way,
\begin{equation}\label{Asympt2}
 |^t(DT_\delta)^{-1}\vec{n}| = 1-(\vec{n},DR.\vec{n})\delta + o(\delta). 
\end{equation}
We also have 
\begin{equation}\label{Asympt3}
 |det~DT_\delta| = det~DT_\delta = 1+\delta (div~R) + o(\delta),
\end{equation}
all the $o(\delta)$ being uniform in $x\in\bar\Omega$. 
Since $R\in C^1(\bar\Omega)$, we get 
$$ Q_\delta(u)=\frac{\displaystyle (1+O(\delta))\int_\Omega (|\nabla u|+|u|\,dx) }
                    { \displaystyle (1+O(\delta))\int_{\p\Omega} |u|\,dH^{N-1}} 
             = (1+O(\delta))Q(u),  $$
as we wanted to prove. Theorem \ref{thm4} then easily follows.

\medskip

We now assume that $\lambda<1$. Since then $\limsup_{\delta\to 0}\lambda_\delta<1$, it follows from Theorem \ref{thm2} that there exists a nonnegative extremal $v_\delta\in BV(\Omega_\delta)$ for $\lambda_\delta$ normalized by $\int_{\p\Omega_\delta}  v_\delta\,dH^{N-1}=1$.
Let $u_\delta = v_\delta\circ T_\delta\in BV(\Omega)$. Then the sequence $(u_\delta)$ is bounded in $BV(\Omega)$. Indeed, according to (\ref{Asympt1}) and  (\ref{Asympt3}), we have 
 \begin{equation*}
 \begin{split} 
   \int_\Omega |\nabla u_\delta| + \int_\Omega u_\delta \,dx 
  & = \int_{\Omega_\delta} |(DT_\delta^{-1})^{-1}\nu_{v_\delta}||det~DT_\delta^{-1}||\nabla v_\delta| 
      + \int_{\Omega_\delta} v_\delta |det~DT_\delta^{-1}| \,dx \\
  & = (1+O(\delta))\int_{\Omega_\delta} |\nabla v_\delta| + v_\delta \,dx 
   = (1+O(\delta))\lambda_\delta \\
  & = (1+o(1))\lambda.
 \end{split}
 \end{equation*} 
 There thus exists a nonnegative $u\in BV(\Omega)$ such that $u_\delta\to u$ in $L^1(\Omega)$. Moreover, as in the proof of theorem \ref{thm2}, 
 \begin{equation*}
 \begin{split} 
  |\nabla u_\delta & |\rightharpoonup \mu\ge |\nabla u| + \sum_{i\in I} \nu_i\delta_{x_i}, \\
  |u_\delta|\,dH^{N-1} & \rightharpoonup \nu=|u|\,dH^{N-1} + \sum_{i\in I} \nu_i\delta_{x_i}. 
 \end{split}
 \end{equation*} 
 We can now obtain
 \begin{equation*}
 \begin{split} 
  \lambda & = \lim_{\delta\to 0} \lambda_\delta = \lim_{\delta\to 0} Q_\delta(v_\delta) = \lim_{\delta\to 0} (1+O(\delta))Q(u_\delta) 
          \ge \frac{\displaystyle\int_\Omega |\nabla u| + \sum_{i\in I} \nu_i+ \int_\Omega u \,dx }
                     {\displaystyle\int_{\p\Omega} u\,dH^{N-1} + \sum_{i\in I} \nu_i} \\
          & \ge \frac{\displaystyle\lambda\int_{\p\Omega} u\,dH^{N-1} + \sum_{i\in I} \nu_i}
                     {\displaystyle\int_{\p\Omega} u\,dH^{N-1} + \sum_{i\in I} \nu_i},
 \end{split}
 \end{equation*} 
 i.e. $\lambda\sum_{i\in I} \nu_i \ge  \sum_{i\in I} \nu_i$. Since $\lambda<1$, we must have $\nu_i=0$ for all $i\in I$, so that 
 $$ 1 = \int_{\p\Omega} v_\delta\,dH^{N-1} = \int_{\p\Omega} u_\delta\,dH^{N-1} + o(1) = \int_{\p\Omega} u\,dH^{N-1} + o(1).  $$
 Using the inferior semi-continuity of the total variation, we can now write 
 \begin{equation*}
 \begin{split} 
  \lambda & = \lim \lambda_\delta = \lim Q_\delta(v_\delta) = \lim~ (1+O(\delta))Q(u_\delta) 
           \ge \frac{\displaystyle\int_\Omega |\nabla u| + \int_\Omega u \,dx }{\displaystyle\int_{\p\Omega} u\,dH^{N-1} } 
           \ge \lambda. 
 \end{split}
 \end{equation*} 
 Hence $u$ is an eigenfunction for $\lambda$ and 
 \begin{equation}\label{StrongConv}
 \begin{split}
  \int_\Omega |\nabla u_\delta| & \to \int_\Omega |\nabla u|, \\
  \int_{\p\Omega} u_\delta\,dH^{N-1} & \to \int_{\p\Omega} u\,dH^{N-1}.
 \end{split}
 \end{equation}
 We now prove the formula for the derivative (\ref{ShapeDer}). We first get using (\ref{Asympt1})-(\ref{Asympt3}) that 
 \begin{equation*}
 \begin{split} 
  Q_\delta(u) 
  & = \frac{\displaystyle \int_\Omega \left(1+\delta f(\nu)+o(\delta)\right)|\nabla u| + \int_\Omega (1+\delta div~R+o(\delta))u \,dx }
         {\displaystyle \int_{\p\Omega} (1+\delta f(\vec{n})+o(\delta))u\, dH^{N-1} } \\
  & = \frac{\displaystyle \lambda + \delta \left(\int_\Omega f(\nu)|\nabla u| + u\,div~R\,dx \right) + o(\delta) }
           {\displaystyle 1+\delta \int_{\p\Omega} f(\vec{n})u\,dH^{N-1} + o(\delta) } \\
  & = \lambda + \delta\left( \int_\Omega (f(\nu)|\nabla u| + u\,div~R\,dx) - \lambda\int_{\p\Omega} f(\vec{n})u\,dH^{N-1}\right) + o(\delta), 
 \end{split}
 \end{equation*} 
 where 
 \begin{equation}\label{def_f}
  f(X)=div~R - (X,DR.X),~X\in\R^n.
 \end{equation}
 Hence
 \begin{equation}\label{Der1}
 \begin{split} 
  \lambda_\delta-\lambda \le & ~Q_\delta(u)-\lambda \\
    = & ~\delta\left( \int_\Omega (f(\nu)|\nabla u| + u\,div~R\,dx) - \lambda\int_{\p\Omega} f(\vec{n})u\,dH^{N-1}\right) + o(\delta). 
 \end{split}
 \end{equation} 
 It remains to prove the opposite inequality. Letting $\nu_\delta\equiv \nu_{u_\delta}$, we obtain, using (\ref{Asympt1}), (\ref{Asympt2}), (\ref{Asympt3}) and the strong convergence $u_\delta\to u$ in $L^1(\Omega)$, that 
 \begin{equation*}
 \begin{split} 
   Q_\delta(u_\delta) 
  & = \frac{\displaystyle\int_\Omega \left\{ 1+\delta f(\nu_\delta)+o(\delta)\right\} |\nabla u_\delta| 
             + \int_\Omega (1+\delta\,div~R + o(\delta))u_\delta \,dx }
             {\displaystyle\int_{\p\Omega} |u_\delta|\,dH^{N-1} +\delta\int_{\p\Omega} f(\vec{n})u_\delta\,dH^{N-1} +o(\delta)} \\
  & = \frac{\displaystyle\int_\Omega (|\nabla u_\delta|+u_\delta\,dx) 
               + \delta \int_\Omega \{f(\nu_\delta)|\nabla u_\delta|+(div~R)u\,dx\} + o(\delta)}
           {\displaystyle\int_{\p\Omega} u_\delta\,dH^{N-1} +\delta\int_{\p\Omega} f(\vec{n})u\,dH^{N-1} +o(\delta)}. 
 \end{split}
 \end{equation*} 
 We can rewrite (\ref{StrongConv}) as 
 \begin{equation}\label{StrongConv2}
  \int_{\bar\Omega} |\nabla\bar u_\delta| \to \int_{\bar\Omega} |\nabla\bar u|, 
 \end{equation}
 where $\bar u_\delta$ (resp. $\bar u$) denotes the extension of $u_\delta$ (resp. $u$) to $\R^n\backslash\bar\Omega$ by $0$. Independently, we clearly
 have the weak  convergence of $\nabla \bar u_\delta$ to $ \nabla \bar u$. We can thus apply Reshetnyak' theorem \cite{Reshetnyak,LuMo,AFP} to get that 
 $$ \int_{\bar\Omega} g(x,\nu_\delta(x)) |\nabla \bar u_\delta| \to  \int_{\bar\Omega} g(x,\nu(x)) |\nabla \bar u| $$
 for any continuous function $g:\bar\Omega\times S \to \R$, where $S$ denotes the unit sphere of $\R^n$. In particular 
 $$ \int_\Omega f(\nu_\delta)|\nabla u_\delta| \to \int_\Omega f(\nu)|\nabla u|. $$
 Hence  
 $$ Q_\delta(u_\delta) = 
   Q(u_\delta) + \delta\left\{ \int_\Omega (f(\nu)|\nabla u| + u\,div~R\,dx) - \lambda\int_{\p\Omega} f(\vec{n})u\,dH^{N-1}\right\} + o(\delta).
 $$
 We now have
 \begin{equation}\label{Der2}
 \begin{split} 
  \lambda_\delta-\lambda 
  \ge & ~Q_\delta(u_\delta)-Q(u_\delta) \\
  = & ~\delta\left( \int_\Omega (f(\nu)|\nabla u| + u\,div~R\,dx) - \lambda\int_{\p\Omega} f(\vec{n})u\,dH^{N-1}\right) + o(\delta).
 \end{split}
 \end{equation} 
 We deduce from (\ref{Der1}) and (\ref{Der2}) and the uniqueness of $u$ that the map $\delta\to \lambda_\delta$ is differentiable at $\delta=0$ with 
 \begin{equation}\label{Der3}
 \begin{split} 
  \lambda_\delta'(0) = \int_\Omega (f(\nu)|\nabla u| + u\,div~R\,dx) - \lambda\int_{\p\Omega} f(\vec{n})u\,dH^{N-1}.
 \end{split}
 \end{equation} 
 As there always exists an eigenset $A\subset\bar\Omega$, i.e. a set of finite perimeter that attains the infimum in (\ref{BestConst_3}), and since $u$
 is by hypothesis the only normalized eigenfunction for $\lambda$, we have $u=|A\cap\p\Omega|^{-1}\chi_A$. It follows from geometric measure theory
 that $|\nabla\chi_A|=|A\cap\p\Omega|^{-1} H^{N-1}_{|\p^*A}$ (see \cite{AFP,EvaGar,Ziemer}). Recalling the definition (\ref{def_f}) of $f$ and using the Green' formula for sets of finite perimeter, we can now rewrite (\ref{Der3}) as (\ref{ShapeDer}). 

\vspace{1cm}

\noindent\textit{\textbf{Acknowledgments.}} The author aknowledges the support of the grants FONCYT PICT 03-13719 (Argentina) and would like to express his gratitude to Prof. J.D. Rossi for its help.

\end{document}